\documentclass[1p, times, twocolumn]{elsarticle}
\usepackage{amssymb}
\usepackage{subfigure}

\usepackage{url}
\usepackage{amsmath, amssymb, color}
\usepackage{graphicx}

\newtheorem{lemma}{Lemma}

\newtheorem{theorem}{Theorem}

\newtheorem{claim}{Claim}

\begin{document}
\begin{frontmatter}
\title{A Simple Policy for Multiple Queues with Size-Independent
Service Times}

\author[umn]{Yuhang Liu}
\ead{liuxx839@umn.edu}

\author[umn]{Zizhuo Wang}
\ead{zwang@umn.edu}

\address[umn]{Industrial and Systems Engineering, University of
Minnesota, MN, 55455}
\begin{abstract}
We consider a service system with two Poisson arrival queues.  A
server chooses which queue to serve at each moment. Once a queue is
served, all the customers will be served within a fixed amount of
time. This model is useful in studying airport shuttling or certain
online computing systems. We propose a simple yet optimal
state-independent policy for this problem which is not only easy to
implement, but also performs very well.
\end{abstract}
\end{frontmatter}

\section{Introduction}
\label{section:introduction} Finding efficient and practical service
plans for multi-queue systems has always been an important question
for both researchers and practitioners. In this paper, we consider a
special type of two-queue system that has the following features: 1)
the customer arrivals in both queues are Poisson processes; 2) there
is only one server; 3) once a queue is chosen to be served, all the
customers in that queue will be served in a fixed amount of time.
The decision in such a system is to decide at each time, which queue
should be served and the overall objective is to minimize the total
expected delays of the customers.

The above model covers a wide range of problems that may occur in
practice. For example, in an airport transportation system, a
shuttle bus picks up passengers from multiple locations (e.g.,
different hotels or rental car locations) and drops them off at the
terminal. Due to the route constraints, the shuttle bus can only go
to one location each time (i.e., a cyclic route is not permissible).
Customers arrive at each location with a certain rate and the
service provider has to decide its service schedule (the schedule
could be either state-dependent or independent) to minimize the
total delays of all customers. Given sufficiently large capacity,
the shuttle can pick up all the passengers waiting in the location
that it serves and the service time does not depend on the job size.
Another example is an online computing service system, in which a
server provides computational services for several types of
customers. In some applications, the service time mainly depends on
the setup time (software initializations and warm up, etc) and
depends little on the job size. In this case, again, the decision
for the server is to decide which type of customer to serve at each
time period, with an overall objective of minimizing the total
delays of jobs.

In this work, we first formulate this problem as a Markov Decision
Process (MDP) in Section \ref{section:model}. Although the MDP
formula is a precise description of the actual problem and the
policy obtained therefore is optimal, the optimal policy has the
feature of being state-dependent. However, in practice, such a
feature might be undesirable due to the limit knowledge of the state
of the system and the implementation constraints. To solve this
problem, in Section \ref{section:heuristics}, we propose a simple
policy that is independent of the current state. Among the
state-independent policies, we prove that the optimal one must
consist of cycles of a fixed schedule, and each cycle consists of
serving the slow-arriving queue once followed by serving the
fast-arriving queue multiple ($k$) times. We then show how to choose
the optimal $k$. We give an explicit formula for the optimal $k$. In
particular, we show that when the discount factor across periods
approaches $1$, $k$ will be approximately $\sqrt{2r}-1$, where $r$
is the ratio of the arrival rates between the fast-arriving queue
and the slow one. This result quantifies the optimal strategy for
the service providers of how often should the slow-arriving queue be
served, and thus is insightful for practitioners. We implement this
strategy in Section \ref{section:numerical} and show that it
performs well.

Before we proceed to the models, we review some related literature
in the following.

\section{Literature Review}
Our work is related to the batched service systems in queueing
theory. For a comprehensive review of batched service system, we
refer the readers to the books by \citet{chaudhry}, \citet{cooper}
and \citet{gross}.

One research question that has attracted a lot of attention in
queueing community is how to optimally schedule a multi-class queue.
In such problems, there is a single server. Multiple classes of
customers are served using this server with different service time
distributions. At each time the server becomes idle, it must choose
which class of queue to serve. And the objective is to find a
service discipline so as to minimize the total expected cost over a
certain horizon. A classical $c\mu$-rule has been proposed for this
model, see, e.g., \citet{baras2,barask,buyukkoc,yao}. In the
$c\mu$-rule, the priority of jobs are ranked by descending order of
the product of $c$ and $\mu$, where $c$ and $\mu$ refer to the cost
and the inverse of the mean service requirement of each customer
class. It has been proved in a vast literature that the $c\mu$-rule
is optimal under a variety of input assumptions, and this model has
been applied to many practical systems, for example, traffic control
systems (\citet{dunne}). For the sake of brevity, we are not going
to review these results in more detail, but refer the readers to the
references above and thereafter. Our problem is similar to this
stream of studies in terms of the decisions being made, however, we
assume that the jobs are served in a batched fashion and within a
fixed amount of time rather than having an i.i.d. service time. This
difference distinguishes our model from this line of research.

There also has been a vast literature in designing systems where a
shuttle runs between two (or more) terminals. The studies of those
problems can be classified by whether the shuttle is finitely
capacitated (\citet{deb2, deb1}) or not (\citet{ignall, weiss1,
weiss2}); whether the control can be exerted at both terminals
(\citet{deb1, deb2, weiss1, lee}) or only can be exerted at a single
terminal (\citet{ignall, weiss2}). The cost usually consists of a
fixed cost, a per passenger costs, and the waiting costs of the
customers. In our work, we assume there are two queues, and the
service provider can serve one of them at each time. Therefore,
although the background shares some similarity, our model is quite
different from
this body of studies.%
%


\section{Model and Analysis}
\label{section:model} We consider a system of two independent queues
without capacity limits. Customers arrive in each of these two
queues according to Poisson processes, with rates $\lambda_1$ and
$\lambda_2$ respectively. One distinguishing feature of our model is
that at each time the server is idle, it could choose one of the two
queues to serve. And once a queue is chosen to be served, all of its
current customers will be served (and depart the system) within a
deterministic time (However, the customers that arrive after the
service has started will not be served during this service cycle).
This model characterizes situations in which the service time is
largely determined by the setup time but not how much service it
needs to provide. Furthermore, we assume the service time for both
queues are the same. Without loss of generality, we assume that the
service time equals to $1$.

We use $Q_i$, $i=1,2$ to denote the $i$th queue. Due to the
assumption of our problem, we can consider this problem in a number
of discrete time periods. The decision in our model is to decide a
queue to serve at each time period. And the objective is to
determine a service rule such that the expected discounted total
waiting time of all the customers is minimized. To mathematically
capture this objective function, we adopt a Markov Decision Process
(MDP) approach. We use $s = (x,y)$ to denote the state when there
are $x$ people waiting in the first queue and $y$ people waiting in
the second queue. We use $V (x,y)$ to denote the optimal expected
discounted waiting time onward (the cost function) when the state is
$(x,y)$. The Bellman equation for this MDP problem can be written as
\begin{eqnarray}\label{dp}
V(x,y)= \lambda + \min \{\gamma \mathbb E [V(Z_1,Z_2+y)]+y, \gamma
\mathbb E [V(Z_1+x,Z_2)]+x \}.
\end{eqnarray}
where $\lambda = (\lambda_1 + \lambda_2)/2$, $Z_1$ and $Z_2$ are
independent Poisson random variables with parameter $\lambda_1$ and
$\lambda_2$ respectively. In (\ref{dp}), the first term is the
expected waiting time for the arriving customers during the current
period. Note that if the arrivals are Poisson processes, given the
number of arrivals in one period, the exact arrival time during this
period is uniformly distributed. Therefore, the expected waiting
time of the new arrivals are $\lambda_i/2$ for each queue $i$. The
second term in (\ref{dp}) corresponds to the two choices available
at the current period, i.e., to serve either the first queue or the
second queue. When the first queue is served, the state will change
to $(Z_1, Z_2 + y)$ and the expected waiting time occurred in this
time period for the customers currently in the queue is $y$.
Similarly, when the second queue is served, the state will change to
$(Z_1 + x, Z_2)$ and the expected waiting time occurred in this time
period for the customers currently in the queue is $x$. $\gamma$ is
the discount factor for the later periods.

The MDP defined in (\ref{dp}) is an infinite state problem. In order
to solve it exactly, we truncate it into a finite state space
problem and apply the normal procedures to solve it (e.g., value
iteration). The solution method is standard and is thus omitted.

\section{State-Independent Policies} \label{section:heuristics}

In Section \ref{section:model}, we established an MDP model for this
problem. Although the MDP problem is easily solvable, the optimal
control is state-dependent, that is, at each moment, the service
provider needs to know the exact queue length of each queue to make
his decision. However, in practice, a state-independent policy might
be desirable. There are two reasons for this. First, it might be
impractical for the service provider to know exactly how many people
are waiting in each queue at each moment. Consider the airport
shuttle example we gave earlier, when dispatching a shuttle to one
location, it is hard to know how many people are waiting in that
location. Second, in real operations, it is very desirable to have a
fixed service schedule. A fixed service schedule not only simplifies
the service provider's task, but also relaxes the customers by
informing them the next service time. In this section, we study
state-independent policies. We first find the optimal
state-independent policy and then show that it performs quite well
in test problems.

\subsection{Optimal structures of the state-independent policy}
\label{subsec:heuristic_structure}

In this section, we prove the structures of the optimal
state-independent policy. In the following discussion, without loss
of generality, we assume $\lambda_1 \le \lambda_2$. We prove the
following theorem.

\begin{theorem}
\label{thm:optimal_heuristic} The optimal state-independent policy
is cyclic after a certain number of periods. In each cycle, $Q_1$ is
served once followed by serving $Q_2$ $k\ge 1$ times.
\end{theorem}

The proof of Theorem \ref{thm:optimal_heuristic} follows from the
following claims.

\begin{claim}\label{claim:1}
The optimal state-independent policy is cyclic after a certain
number of periods. In each cycle, $Q_1$ is served $k_1$ times
followed by serving $Q_2$ $k_2$ times.
\end{claim}

\begin{claim}\label{claim:2}
$k_1 = 1$.
\end{claim}

{\bf Proof of Claim \ref{claim:1}.} We consider any
state-independent policy $\pi$. We call $\pi$ ``reasonable'' if it
serves each queue infinite many times. It is obvious that the
optimal policy must be a ``reasonable'' one. Now we only focus on
``reasonable'' policies. Consider the following piece of service
schedule $\pi(l_1,l_2)$ ($l_1\ge 0$, $l_2\ge 1$): Serve $Q_1$ $l_1$
times, then serve $Q_2$ $l_2$ times, then serve $Q_1$. Note that any
``reasonable'' policy can be decomposed into such pieces based on
each occurrence of $(Q_2, Q_1)$ in its schedule. Also note that for
fixed $l_1$ and $l_2$, the expected costs during this piece
discounted back to its beginning is a constant. To see this, we
first note that since the prior piece always ends up with serving
$Q_2$ then $Q_1$ (by definition), at the beginning of each piece,
the initial state must have a distribution of $(Z_1,
Z_2^{1}+Z_2^{2})$, where $Z_1$ is a Poisson random variable with
parameter $\lambda_1$ and $Z_2^{i}$ ($i=1,2$) are independent
Poisson random variables with parameter $\lambda_2$. Thus given
$(l_1,l_2)$ fixed, the expected cost of a piece $(l_1, l_2)$ is a
constant. Denote $t_0$ to be the first time that $Q_2$ is served at
$t_0-2$ and $Q_1$ is served at $t_0-1$. Denote $V$ to be the optimal
expected cost from $t_0$. Now if in the optimal policy, the first
piece is $\pi(l_1, l_2)$ (as argued before, it must complete such a
piece at a certain point), then we must have
\begin{eqnarray*}
V = V(l_1,l_2) +  \gamma^{l_1+l_2 + 1} V
\end{eqnarray*}
where $V(l_1, l_2)$ is the expected cost during this piece
discounted back to its beginning. This is because the state after
$l_1+l_2+1$ periods after $t_0$ have exactly the same distribution
as that at time $t_0$ and thus the optimal expected cost onward are
the same. Noting that repeatedly using $\pi(l_1, l_2)$ will achieve
$V$. Therefore, repeatedly using $\pi(l_1,l_2)$ is an optimal policy
from time $t_0$.

Finally, if we put $\pi(l_1, l_2)$ together and break it by each
time we switch from $Q_2$ to $Q_1$, we get the optimal
state-independent policy must be cyclic from $t_0-1$ and in each
cycle, $Q_1$ is served $k_1$ times followed by $Q_2$ served $k_2$
times. \hfill$\Box$

{\bf Proof of Claim \ref{claim:2}.} We prove by showing that for any
cyclic policy with $k_1 > 1$, we can always find another policy that
reduces the expected cost. We consider two cases:

\begin{itemize}
\item Case 1: $k_1 > 2$. Then we compare the following two policies:

\begin{enumerate}
\item Serve $Q_1$ $k_1$ times, followed by serving $Q_2$ $k_2$
times, i.e.. $\begin{array}{cc}
\underbrace{Q_1...Q_1}\underbrace{Q_2...Q_2}\\
k_1 \quad\mbox{  }\quad k_2
\end{array}$
\item Serve $Q_1$ once, then $Q_2$ once, then $Q_1$ $k_1 - 2$ times,
then $Q_2$ $k_2$ times, i.e., $
\begin{array}{cc}
Q_1Q_2\underbrace{Q_1...Q_1}\underbrace{Q_2...Q_2}\\
\mbox{ }\quad\mbox{     }k_1 - 2 \mbox{     }\mbox{  }\quad \mbox{ }
k_2
\end{array}$
\end{enumerate}

Since the time periods covered by policy $1$ and $2$ are both
$k_1+k_2$, therefore to compare the expected cost, we can simply
compare the expected cost within each cycle. The expected cost for
each cycle of the first policy is:
\begin{eqnarray*}
C_1 = \lambda\sum_{i=0}^{k_1+k_2 - 1}\gamma^i + \lambda_2
\sum_{i=1}^{k_1} i\gamma^{i-1} +\lambda_1\gamma^{k_1}
\sum_{i=1}^{k_2}i\gamma^{i-1},
\end{eqnarray*}
where the first term is the waiting time incurred by the new
arrivals in each period, the second term is the waiting time of
$Q_2$ in the first $k_1$ periods, and the third term is the waiting
time of $Q_1$ in the last $k_2$ periods. Similarly, the cost for
each cycle of the second policy is:
\begin{eqnarray*}
C_2 = \lambda\sum_{i=0}^{k_1 + k_2 -1} \gamma^i + \lambda_2 +
\lambda_1\gamma + \lambda_2 \sum_{i=3}^{k_1} (i-2)\gamma^{i-1}
+\lambda_1\gamma^{k_1}  \sum_{i=1}^{k_2}i\gamma^{i-1}.
\end{eqnarray*}
Now we compare $C_1$ and $C_2$, we have $C_1 - C_2 = \lambda_2
\sum_{i=2}^{k_1} 2\gamma^{i-1} -\lambda_1 \gamma > 0
$
when $\lambda_2 > \lambda_1$. Therefore, we proved that any cycle
with $k_1 > 2$ can't be optimal.
\item Case 2: $k_1 = 2$. We consider two further cases:
\begin{itemize}
\item Case 2.1: $k_2 = 1$, then we compare the cycles $Q_1Q_1Q_2$
with $Q_1Q_2Q_2$. Again, since both cycles cover $3$ time periods,
we only need to compare the expected cost within the cycle. In the
cycle $Q_1Q_1Q_2$, the expected cost is $\lambda\sum_{i=0}^2\gamma^i
+ \lambda_2 + 2\lambda_2\gamma +\lambda_1\gamma^2$, whereas in the
cycle $Q_1Q_2Q_2$, the expected cost is $\lambda\sum_{i=0}^2\gamma^i
+ \lambda_2 + \lambda_1\gamma + 2\lambda_1\gamma^2$. It is easy to
see that when $\lambda_2 > \lambda_1$ and $\gamma < 1$, the expected
cost in the first cycle is larger. Therefore, it can't be the
optimal policy.
\item Case 2.2: $k_2\ge 2$. Then we compare the following two
policies:
\begin{enumerate}
\item Serve $Q_1$ twice followed by $Q_2$ $k_2$ times
\item Serve $Q_1$ once, then $Q_2$ once, then $Q_1$ once and then
$Q_2$ $k_2 - 1$ times
\end{enumerate}
Again, we only need to compare the cost within these two cycles. In
the first cycle, the expected cost is:
\begin{eqnarray*}
C_1 = \lambda\sum_{i=0}^{k_2 + 1} \gamma^i + \lambda_2 +
2\lambda_2\gamma + \lambda_1 \gamma^2 \sum_{i=1}^{k_2}i\gamma^{i-1}
\end{eqnarray*}
where the first term is the waiting time due to the arrivals within
each time period, the second and third term are the waiting time in
$Q_2$ in the first and second time periods and the last term is the
waiting time in $Q_1$ in the last $k_2$ time periods. Similarly, the
expected cost in the second cycle is:
\begin{eqnarray*}
C_2 = \lambda \sum_{i=0}^{k_2 +1 } \gamma^i + \lambda_2 +
\lambda_1\gamma + \lambda_2 \gamma^2 +\lambda_1 \gamma^3
\sum_{i=1}^{k_2-1}i\gamma^{i-1}
\end{eqnarray*}
By taking the difference, we have $C_1 - C_2 = 2\lambda_2 \gamma -
\lambda_1\gamma -\lambda_2\gamma^2 +\lambda_1 \sum_{i=1}^{k_2}
\gamma^{i+1} > 0$ when $\lambda_2 > \lambda_1$ and $\gamma < 1$.
Therefore $k_1 = 2$ can't be optimal too and Claim \ref{claim:2} is
proved. \hfill $\Box$
\end{itemize}
\end{itemize}

\subsection{Finding optimal $k$} \label{subsec:find k}

We first compute the expected discounted cost in each cycle.

\begin{lemma}
Consider a cycle of $(Q_1,Q_2,...,Q_2)$ with length $k+1$ (thus $k$
$Q_2$s). The expected discounted cost in this cycle is
\begin{eqnarray}\label{cost in one cycle}
\lambda \sum_{i=0}^k \gamma^i + \lambda_2 + \lambda_1
\sum_{i=1}^{k_2} i\gamma^i.
\end{eqnarray}
\end{lemma}

In (\ref{cost in one cycle}), the first term is the waiting time due
to the arrivals within each period, the second term is the waiting
time in $Q_2$ in the first time period and the last term is the
waiting time in $Q_1$ in the last $k$ periods. The lemma just
follows by summing these costs together.
%
%
%

Now we consider the total expected costs throughout the time
horizon, if a cyclic policy with length $k+1$ is used, each with
cost $C$, then the total discounted cost will be
$C/(1-\gamma^{k+1})$. Therefore we have the following theorem.

\begin{theorem}\label{thm:total_cost} Assuming the initial state is
$(M,\lambda_2)$ where $M >\lambda_2$. Then the total expected cost
of using cyclic policy $(Q_1,Q_2,...,Q_2)$ with length $k+1$ (thus
$k$ $Q_2$s) is
\begin{eqnarray*}\label{total_cost}
C(k)= \frac{\lambda \sum_{i=0}^k \gamma^i + \lambda_2 + \lambda_1
\sum_{i=1}^{k} i\gamma^i}{1-\gamma^{k+1}}.
\end{eqnarray*}
\end{theorem}

{\noindent\bf Remark.} The reason we involve the big $M$ is to make
sure that it is better to serve $Q_1$ in the first period, and the
expected cost in the very first period is the same as in the first
period of each cycle later. In practice, if one focuses on the long
run, then the costs in the first period does not matter too much and
the result still approximately holds, otherwise, a modification is
needed. \hfill$\Box$

Now our task is simply to find $k$ to minimize $C(k)$. The following
theorem shows that $C(k)$ is unimodal in $k$. For the convenience of
notation, we denote $r= \frac{\lambda_2}{\lambda_1} > 1$ in the
following discussions.

\begin{theorem}\label{thm:find_optimal_k}
There exists a unique $k^*$ such that $C(k^*-1) \ge C(k^*)$
 and $C(k^*+1) \ge C(k^*)$.
\end{theorem}
{\noindent\bf Proof.} We compare $C(k)$ and $C(k+1)$. Note that the
$\lambda$ part in $C(k)$ are the same for all $k$ (all equal to
$\frac{\lambda}{1-\gamma}$). Therefore, it suffices to compare
\begin{eqnarray*}
C'(k)  = \frac{\lambda_2 + \lambda_1 \sum_{i=1}^{k}
i\gamma^i}{1-\gamma^{k+1}}
\end{eqnarray*}
with $C'(k+1)$. Further by factoring out $\lambda_1$, we only need
to compare
\begin{eqnarray*}
\bar{C}(k) = \frac{r + \sum_{i=1}^{k} i\gamma^i}{1-\gamma^{k+1}}
\end{eqnarray*}
with $\bar{C}(k+1)$. We have $\bar{C}(k+1) \ge \bar{C}(k)$ if and
only if
\begin{eqnarray}\label{condition1}
(r + \sum_{i=1}^{k}i\gamma^i)(1-\gamma^{k+2}) - (r +
\sum_{i=1}^{k+1}i\gamma^i)(1-\gamma^{k+1}) \le 0.
\end{eqnarray}
By expanding and regrouping terms, we get that (\ref{condition1}) is
equivalent as
\begin{eqnarray}\label{condition2}
r \le (k+1)\sum_{i=0}^k \gamma^i - \sum_{i=1}^k i\gamma^i =
\sum_{i=0}^k(k+1-i)\gamma^i.
\end{eqnarray}
Note that the right hand side of (\ref{condition2}) is increasing in
$k$. And when $k = 0$, the right hand side is $1$ which is less than
$r$, and when $k$ goes to infinity, the right hand also goes to
infinity (since it is greater than $k$). Therefore, there exists a
unqiue integer $k^*$ such that
\begin{eqnarray}\label{condition3}
 \sum_{i=0}^{k^*}(k^*-i)\gamma^i \le r <
 \sum_{i=0}^{k^*+1}(k^*+1-i)\gamma^i.
\end{eqnarray}
Such $k^*$ will satisfy the property stated in the theorem. \hfill
$\Box$

{\noindent\bf Remark.} By Theorem \ref{thm:find_optimal_k}, finding
the optimal $k$ in the state-independent policy simply becomes to
find $k$ that satisfies (\ref{condition3}). This can be done very
fast. Moreover, we can see that $k$ only depends on $\lambda_1$ and
$\lambda_2$ through $r = \frac{\lambda_2}{\lambda_1}$. When $r$ is
larger, $k$ is also larger, and vice versa. This is consistent with
our intuition that if the arrival rate of $Q_2$ is much larger than
that of $Q_1$, we should serve it more frequently.

\begin{figure}[htbp]
\centering\subfigure[Optimal $k^*$]{
\includegraphics[width=2.2in]{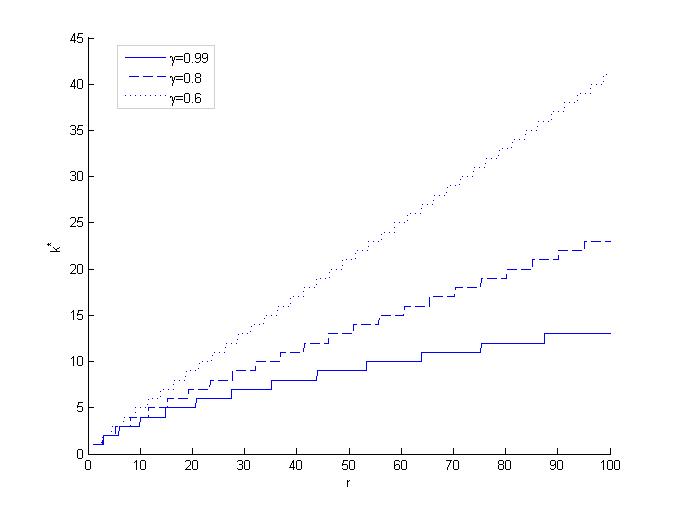}
\label{fig:optimal_k} } \subfigure[$C(k)$ with $\lambda_1=1$,
$\lambda_2=9$, $\gamma=0.8$]{
\includegraphics[width=2.2in]{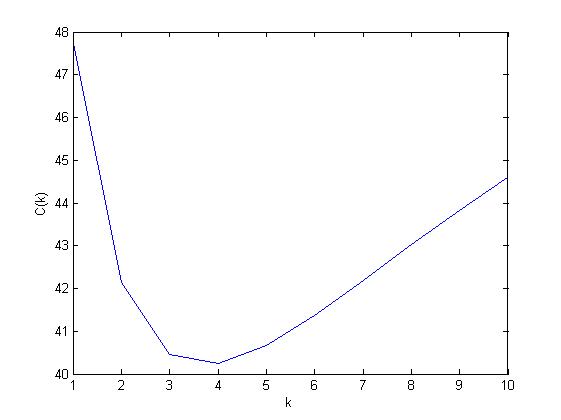}
\label{fig:c} } \caption{Illustrations of $k^*$ and $C(k)$}
\label{fig:illustration}
\end{figure}

In Figure \ref{fig:illustration}, we show an example of the optimal
$k^*$ and the function $C(k)$.

The next theorem points out the relationship between the optimal $k$
and the discount factor $\gamma$ and two important limit cases. The
proof follows immediately from (\ref{condition3}).
\begin{theorem}\label{thm:limitcase}
The optimal $k$ is decreasing with the discount factor $\gamma$.
Furthermore, when $\gamma\rightarrow 1$, $k\sim\sqrt{2r}-1$ and when
$\gamma \rightarrow 0$, $k\sim r$.
\end{theorem}
Theorem \ref{thm:limitcase} is quite instrumental. It directly links
the discount factor to the optimal frequency of serving each queue.
In particular, when discount factor is small, then we are more
focused on the current state, and we will choose the longer queue to
serve at each moment. Under the state-independent assumption, this
means that if the ratio between the arrival rates is $r$, then the
ratio of the frequency of serving the two queues should be about $r$
(i.e., serve the slow-arrival queue only when it accumulates the
same amount of customers as the long queue). On the contrary, when
the discount factor goes to $1$, we have to be careful about the
accumulation effect, that is, the accumulated waiting time for $k$
period for one queue is roughly of order $k^2$, and we want to serve
the slow-arrival queue if $k^2$ is of order $r$. Therefore, the
ratio of the frequency between serving the two queues should be
about $O\left(\sqrt{r}\right)$. This rule, although simple, may
provide important guidance for practitioners when they decide how
frequent to serve each queue.

\section{Numerical experiments}
\label{section:numerical}

In this section, we conduct numerical experiments to verify our
findings. The results are shown in Table \ref{table:result}.

\begin{table}[htbp]
\centering
\begin{tabular}{|c|c|c|c|c|c|c|c|c|c|}
\hline $\gamma$  & $r$ & $k^*$ &  $C(1)$ & $C(r)$ & $C(k^*)$ & OPT &
Gap(1) &Gap($r$) & Gap($k^*$) \\ \hline
0.6 &  1 &1& 5.00& 5.00&5.00& 4.62&8.29\%  &8.29\%  &8.29\%   \\
0.6 &  3 &2& 10.63& 10.71&10.51& 9.93&6.98\%  &7.81\%  &5.82\%   \\
0.6 &   5 &3& 16.25& 15.76& 15.51& 14.91&8.96\%  &5.68\%  &3.97\%   \\
0.6 &  9 &4& 27.50& 25.15&24.95& 24.51& 12.20\%  &2.63\%  &1.82\%   \\
\hline
0.8 & 1 &1& 10.00& 10.00&10.00& 8.85&13.04\%  &13.04\%  &13.04\%   \\
0.8 &  3 &2& 20.56& 21.21&20.41& 18.47&11.28\%  &14.80\%  &10.49\%   \\
0.8 &  5&2& 31.11& 31.12&29.51& 27.27&14.06\%  &14.08\%  &8.18\%   \\
0.8 & 9&4& 52.22& 49.07&46.20& 43.93&18.86\%  &11.68\%  &5.16\%   \\
\hline
0.99 & 1 &1& 200.00& 200.00&200.00&167.86&19.14\%  &19.14\%  &19.14\%   \\
0.99 & 3 &2& 399.50& 424.12&399.33& 342.91&16.35\%  &23.43\%  &16.30\%   \\
0.99 & 5&2& 599.00& 630.82&565.33& 500.30&19.49\% &  25.75\% &12.86\% \\
0.99 & 9&3& 997.99& 1032.23&871.87& 718.22&25.37\% &  29.60\% &
9.75\% \\ \hline
\end{tabular}
\caption{Numerical results of our policies} \label{table:result}
\end{table}

In Table \ref{table:result}, the first and second column show the
$\gamma$ and $r$ used. The third column is the optimal $k^*$
computed from (\ref{condition3}). The fourth to sixth column are the
expected costs when one uses a state-independent policy with serving
$Q_1$ once followed by serving $Q_2$ $1$, $r$ or $k^*$ time
respectively. The reason we choose $k=1$ and $k=r$ to compare is
they are the simplest choices thus might be chosen in practice. The
next column OPT is the optimal expected cost in the MDP. Note that
OPT provides a lower bound of the cost for all state-independent
policies. The last three columns in Table \ref{table:result} show
the cost gap between the state-independent policies and $OPT$. For
all the numbers, we use $(M, \lambda_2)$ as our initial state, which
conforms to Theorem \ref{thm:total_cost}.

In Table \ref{table:result}, we see that the policy $k^*$ indeed
performs much better than other choices. In particular, we see that
choosing $k=1$ and $k=r$ usually result in comparable gaps from
$OPT$, while choosing $k^*$ results in about half of the gaps.
Furthermore, the gap becomes smaller when $r$ becomes larger, since
if one queue has much faster arrival than the other, in either
state-independent or dependent policy, one has to serve that queue
more frequently, and the discrepancies between these two policies
are smaller. Also, the gap becomes larger when $\gamma$ becomes
larger, potentially due to that if the performance is evaluated in a
long run, the ability to adjust to the current state is more
important.
%


\bibliographystyle{elsarticle-num-names}
\bibliography{thesis_ref}

\end{document}